\documentclass[number,citesort,MSNbibl,dvips]{arxbj}
\usepackage{upgreek}
\usepackage{graphicx}

\aid{0}
\volume{18}
\issue{1}
\pubyear{2012}
\firstpage{100}
\lastpage{118}
\doi{10.3150/10-BEJ339}

\makeatletter
\newcommand{\eqref}[1]{(\ref{#1})}
\newcommand{\xrightarrow}[1]{\stackrel{#1}{\rightarrow}}

\newcommand{\E}{{\mathsf{E}}}
\newcommand{\N}{{\mathbb{N}}}
\renewcommand{\P}{{\mathsf{P}}}
\newcommand{\R}{{\mathbb{R}}}
\renewcommand{\epsilon}{{\varepsilon}}
\newcommand{\cov}{\operatorname{\mathsf{cov}}}
\newcommand{\corr}{\operatorname{\mathsf{corr}}}
\newcommand{\ind}{\mathbb{I}}
\newcommand{\dist}{\operatorname{\mathsf{dist}}}

\newtheorem{theorem}{Theorem}
\newtheorem{proposition}{Proposition}
\newtheorem{lemma}{Lemma}
\newproclaim{definition}{Definition}
\newtheorem{corollary}{Corollary}
\makeatother

\begin{document}
\begin{frontmatter}

\title{Central limit theorems for the excursion set volumes of weakly
dependent random fields}
\runtitle{CLT for excursion sets of random fields}

\begin{aug}
\author[1]{\inits{A.}\fnms{Alexander}~\snm{Bulinski}\thanksref{1}\ead[label=e1]{bulinski@mech.math.msu.su}},
\author[2]{\inits{E.}\fnms{Evgeny}~\snm{Spodarev}\thanksref{2,e2}\ead[label=e2,mark]{evgeny.spodarev@uni-ulm.de}}\break
\and~%
\author[2]{\inits{F.}\fnms{Florian}~\snm{Timmermann}\corref{}\thanksref{2,e3}\ead[label=e3,mark]{florian.timmermann@uni-ulm.de}}
\runauthor{A. Bulinski,
E. Spodarev
and
F. Timmermann}
\address[1]{Moscow State University, Department of
Mathematics and Mechanics, 119991 Moscow, Russia.\\\printead{e1}}
\address[2]{Ulm University, Institute of Stochastics, 89069 Ulm,
Germany.\\
\printead{e2};
\printead*{e3}}
\end{aug}

\received{\smonth{5} \syear{2010}}
\revised{\smonth{8} \syear{2010}}

%
\begin{abstract}
The multivariate central limit theorems (CLT) for the volumes of excursion
sets of stationary quasi-associated random fields on
$\mathbb{R}^d$ are proved. Special attention is paid to Gaussian
and shot noise fields. Formulae for the covariance matrix of the
limiting distribution are provided. A~statistical version of the
CLT is considered as well. Some numerical results are also discussed.
\end{abstract}

%
\begin{keyword}
\kwd{central limit theorem}
\kwd{dependence conditions}
\kwd{excursion sets}
\kwd{random fields}
\end{keyword}

\end{frontmatter}
%

\section{Introduction}
An important research domain of modern probability theory is the
investigation of geometric characteristics of random surfaces (see, e.g.,
\cite{Adler1,Adler2,Belyaev}).
The origin of interest often roots not only in pure mathematical
challenges but
also in various applications, including those in industry. We
mention one motivating example for our study.

The contemporary method of papermaking goes back to the Han Dynasty
period. Nowadays,
the method is essentially the same, but machines
in modern
pulp and paper mills operate much faster. The
surface structure of the paper during the forming process
determines the quality of the production.

To model the paper surface, stationary random fields, say, shot
noise (cf. \cite{Brown}) or Gaussian, can be a reasonable first
choice. Comparing by eye real paper image data and simulated
realizations of such fields, one easily concludes that the
similarities are striking. But it is hard to quantify how
different these two images really are. To test whether the
available image data originate from a realization of a specified
stationary random field, the excursion sets can be considered.


We prove the central limit theorem (CLT) for volumes of excursion sets
of a stationary
field  $X=\{X(t), t\in T\}$,  $T\subset\mathbb{R}^d$, to
characterize the surface generated by $X$. It is reasonable to
assume that the field $X$ could possess a dependence structure
more general than positive or negative association used in a
number of stochastic models; see, for example,~\cite{Bulinski}. Our main
results yield uni- and multivariate CLT for quasi-associated
random fields. The CLT is generalized in \cite{Ivanov}, page 80,
having been obtained by other methods for volumes of excursion sets of
stationary and isotropic Gaussian random fields. We also discuss
the consistent estimators for the asymptotic covariance matrix that arises
in the limiting distribution.

Note that we do not tackle here the interesting problems concerning
the study of moving levels for excursion sets,
the estimate of the convergence rate to the limit law
and the analysis of the functionals in Gaussian random fields
based on the Dobrushin--Major techniques.
In this regard, we refer, to \cite{Ivanov,Leonenko1,Leonenko2,Leonenko3}.



As to the problem of characterizing the paper quality taking into account
the ``hills'' and ``valleys'' of its surface
discernible with the help of microscope, it is by no means simple.
In fact, we have to specify the
admissible (average) number of such hills along with their size.
Moreover, the thickness of the paper should be controlled as well
(no holes or high peaks). Thus the study of
the excursion sets for random fields is the first natural step to
investigate such random
surfaces. The application to paper surface image data will
appear in a separate paper.

The present paper is organized as follows. Section
\ref{sectPreliminary} provides preliminaries on dependence
concepts related to association and excursion sets for random
fields. The CLT for the volumes of excursions of quasi-associated
stationary random fields over one or finitely many levels are
formulated and proved in Section \ref{sectCLT}. The special cases
of stationary shot noise and Gaussian random fields are treated in
more detail. Section \ref{sectStatCLT} contains a statistical
version of the limit theorems mentioned above where the (unknown)
limiting covariance matrix is consistently estimated. Numerical results
illustrating the limit theorem of Section \ref{sectCLT} are given in
Section \ref{discussion}. Finally, we conclude with the discussion of
some open problems.

\vspace*{2pt}\section{Preliminaries}\label{sectPreliminary}\vspace*{2pt}

In this section, we recall some dependence concepts for systems of
random variables. Various examples can be found in
\cite{Bulinski}. After that, we introduce the excursion sets that
are the main objects of this study. Then we consider the sequences
of regular growing sets forming observation windows.

\vspace*{2pt}\subsection{Dependence concepts for random fields}\vspace*{2pt}

Consider a family, $X=\{X(t), t\in T\},$ of real-valued random
variables, $X(t),$ defined on a probability space,
$(\Omega,\mathcal{F},\mathsf{P})$. A set, $T,$ will be a subset of
$\mathbb{R}^d$ or $\mathbb{Z}^d$. For $I\subset T$ let
$X_I=\{X(t), t\in I\}$. Introduce the class $ \mathcal{M}(n)$
consisting of real-valued, bounded, coordinate-wise non-decreasing
Borel functions on $ \mathbb{R}^n, n \in\mathbb{N} $. The
cardinality of a finite set, $U,$ will be denoted by $\operatorname{\mathsf{card}}
U$.

\begin{definition}\label{PANA}
A real-valued random field $X=\{X(t), t\in T\}$ is called
\textit{positively associated} $($we write $X\in\mathsf{PA})$ if,
for
every disjoint finite set $ I,J \subset T$ and any functions $f\in
\mathcal{M}(\operatorname{\mathsf{card}} I)$ and $g \in\mathcal{M}(\operatorname{\mathsf{card}} J),
$ one has
%
\begin{equation}\label{PA}
\operatorname{\mathsf{cov}} (f(X_I),g(X_J) ) \geq0.\vadjust{\goodbreak}
\end{equation}
\end{definition}

Here, we use any permutation of (coordinates of) the column vector
$(X(t_1),\ldots ,X(t_n))^\top$ for $X_I$,
$I=\{t_1,\ldots,t_n\}\subset T$ (and the analogous notation is
employed for $X_J$); $\top$ stands for transposition.
Definition~\ref{PA}, given for any (not necessarily disjoint)
subsets $I$ and $J\subset T$, introduces the family of \textit{associated} random variables ($X\in\mathsf{A}$). The change of the
sign of inequality \eqref{PA} leads to the definition of {negative
association} (one writes $X\in\mathsf{NA})$. Clearly, $X\in\mathsf{A}$
implies $X\in\mathsf{PA}$. Any collection of independent random
variables is automatically $\mathsf{PA}$ and $\mathsf{NA}$. Due to Pitt
\cite{Pitt}, a Gaussian family $X=\{X(t), t\in T\}$ of random
variables is associated if and only if $\operatorname{\mathsf{cov}}(X(s),X(t))\geq
0$ for all $s,t\in T$. For such families, the concepts of $\mathsf{A}$
and $\mathsf{PA}$ coincide. A~theorem by Joag-Dev and Proschan
\cite{JDP} states that a Gaussian family $X=\{X(t), t\in T\}\in
\mathsf{NA}$ if and only if $\operatorname{\mathsf{cov}}(X(s),X(t))\leq0$ for $s,t\in
T$, $s\neq t$.\looseness=-1

Let $\mathit{BL}(n)$ denote the class of \textit{bounded Lipschitz} functions
$f\dvtx \mathbb{R}^n\to\mathbb{R}$ ($n\in\mathbb{N}$) and
\[
\operatorname{\mathsf{Lip}} (f) = \sup_{x\neq y}
\frac{|f(x)-f(y)|}{ \|x-y \|_1}<\infty, \qquad
 \|x \|_1 = \sum_{k=1}^n|x_k|,
\ x=(x_1,\ldots,x_n)^\top\in\mathbb{R}^n.
\]
Since all norms are
equivalent in $\mathbb{R}^n$, we sometimes use the Euclidean norm
$ \|x \|_2 =   (\sum_{k=1}^n x^2_k )^{1/2}$ and
the supremum norm $ \|x \|_\infty= \max_{k=1,\ldots,n}
|x_k| $ of $x\in\mathbb{R}^n$ for the sake of convenience.

\begin{definition}\label{def:QA}
A random field $X=\{X(t), t\in T\}$ consisting of random
variables $X(t)$ with $\mathsf{E}X(t)^2<\infty$ is called
\textit{quasi-associated} $(X\in\mathsf{QA})$ if
%
\begin{equation}\label{QA}
 | \operatorname{\mathsf{cov}}(f(X_I),g(X_J)) | \leq\operatorname{\mathsf{Lip}}
 (f) \operatorname{\mathsf{Lip}}  (g)\sum_{s\in I} \sum_{t\in J} | \operatorname{\mathsf{cov}}(X(s),X(t)) |
\end{equation}
for all disjoint finite sets $I,J\subset T$ and any Lipschitz
functions
\[
f\dvtx \mathbb{R}^{\operatorname{\mathsf{card}}  I}\to\mathbb{R}   \quad
\mbox{and}\quad
 g\dvtx \mathbb{R}^{\operatorname{\mathsf{card}}  J} \to\mathbb{R}.
 \]
\end{definition}

If $X\in\mathsf{PA}$ or $X\in\mathsf{NA}$ and $\mathsf{E}X(t)^2<\infty$
for all $t\in T,$ then \eqref{QA} holds as was proved in~\cite{BS}.
Every Gaussian random field $X$ (with covariance function taking
both positive and negative values) is quasi-associated; see
\cite{Shash} and references therein.

\begin{definition}\label{BL(theta)}
A real-valued random field $ X=\{X(t), t\in\mathbb{Z}^d\} $ is
called \textit{$(\mathit{BL},\theta)$-dependent}  $(X\in(\mathit{BL},\theta))$ if
there exists a non-increasing sequence $ \theta=
(\theta_r)_{r\in\N},  \theta_r \rightarrow0 $ as $ r
\to\infty$, such that, for any finite disjoint sets $I$,
$J\subset\mathbb{Z}^d $ with $ \dist(I,J) = r$ and
any functions $ f\in \mathit{BL}(\operatorname{\mathsf{card}}  I)$, $g\in \mathit{BL}(\operatorname{\mathsf{card}}
J),$ one has
%
\begin{equation}\label{BL}
 | \operatorname{\mathsf{cov}}(f(X_I),g(X_J)) | \leq\operatorname{\mathsf{Lip}} (f) \operatorname{\mathsf{Lip}} (g)(\operatorname{\mathsf{card}}  I \wedge\operatorname{\mathsf{card}}  J)  \theta_r,
\end{equation}
where $\dist(I,J) = \min\{ \|s-t \|_{\infty}\dvt s\in
I, t\in J\}$.\vadjust{\goodbreak}
\end{definition}

If $X=\{X(t), t\in\mathbb{Z}^d\}\in\mathsf{QA}$, then
$X\in(\mathit{BL},\theta)$ whenever the
Cox--Grimmett coefficient
%
\begin{equation}\label{CG}
u_r:= \sup_{s\in\mathbb{Z}^d} \sum_{t\dvt \|s-t\|_{\infty} \geq
r}|\operatorname{\mathsf{cov}}(X_s,X_t)|
\end{equation}
tends to zero as $r\to\infty$. In this case, one can take $\theta
_r=u_r$ in \eqref{BL}.

For a random field $X=\{X(t), t\in\mathbb{R}^{d}\}$ we use (see
\cite{Bulins}) the following extension of \eqref{BL}. Let
$T(\Delta)=\{(j_1/\Delta,\ldots,j_d/\Delta)\dvt  (j_1,\ldots,j_d)
\in\mathbb{Z}^d\},$ where $\Delta>0$.

\begin{definition}\label{genBL}
A real-valued random field $ X=\{X(t), t\in\mathbb{R}^d\} $ is
called \textit{$(\mathit{BL},\theta)$-dependent} if there exists a
non-increasing function $ \theta= (\theta_r)_{r>0},
\theta_r \rightarrow0 $ as $ r \rightarrow\infty$, such that,
for all~$\Delta$ large enough and any finite disjoint sets $I$,
$J\subset T(\Delta) $ with $ \dist(I,J) = r $, and
any functions $ f\in \mathit{BL}(\operatorname{\mathsf{card}}  I)$, $g\in \mathit{BL}(\operatorname{\mathsf{card}}
J)$, one has
%
\begin{equation}\label{gBL}
 | \operatorname{\mathsf{cov}}(f(X_I),g(X_J)) | \leq\operatorname{\mathsf{Lip}} (f) \operatorname{\mathsf{Lip}} (g)(\operatorname{\mathsf{card}}  I \wedge\operatorname{\mathsf{card}}  J) \Delta^d \theta_r.
\end{equation}
\end{definition}

In many cases, one can use the integral analog of \eqref{CG} for
$\theta_r$. Thus for a (wide-sense) stationary random field
$X=\{X(t), t\in\mathbb{R}^d\}\in\mathsf{QA}$, having
covariance function $R(t)$, \mbox{$t\in\mathbb{R}^d$}, absolutely
directly integrable in the Riemann sense{}
(i.e., when $d=1;$ see, e.g., Feller~\cite{Feller}, page 362.
For $d>1,$ the definition is quite similar. One takes the partition of
$\mathbb{R}^d$ generated by partitions of each coordinate
axis and forms the corresponding upper and lower Riemann sums.),
relation \eqref{gBL}
holds with
%
\begin{equation}\label{int}
\theta_r= 2\int_{\|x\|_\infty\geq r}|R(t)|\,\mathrm{d}t,  \qquad   r>0;
\end{equation}
see \cite{Bulins}. We shall also write
$\theta(X)=\theta_r(X)$ to emphasize that
$\theta$ in \eqref{BL} or \eqref{gBL} refers to the field $X$.

\subsection{Excursion sets}

\begin{figure}

\includegraphics{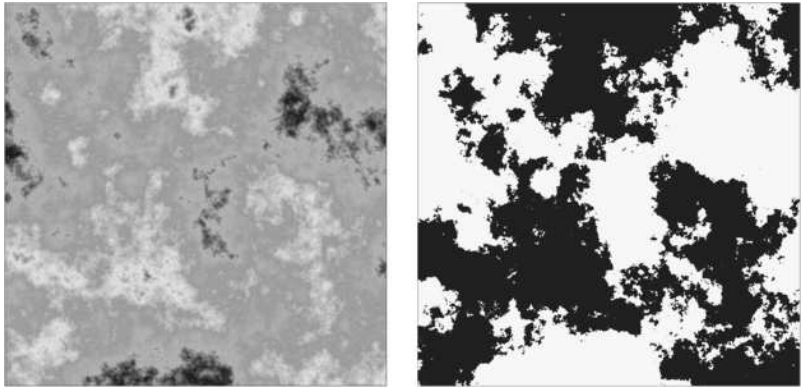}%
\vspace*{-3pt}
\caption{Realization of a stationary centered Gaussian random field $X$
with covariance function $\cov(X(0),X(t))=\exp(- \|t \|
_2)$ (left figure),
bright colours indicate high values of $X$. The excursion set
$A_u(X,T)$ for $u=0$
is shown in black (right figure).}
\label{fig:niveaumenge}\vspace*{-4pt}
\end{figure}

Now we recall the definition of an excursion set and illustrate it
by Figure \ref{fig:niveaumenge}.
%

For a real-valued random field $X=\{X(t), t\in\mathbb{R}^d\},$ we
assume the measurability of~$X(\cdot)$ as a function on
$\mathbb{R}^d\times\Omega$ endowed with the $\sigma$-algebra
$\mathcal{B}(\mathbb{R}^d)\otimes\mathcal{F}$.
\begin{definition}\label{excur}
Let $X$ be a~measurable real-valued function on $\mathbb{R}^d$
and
$T\subset\mathbb{R}^d$ be a~$($Lebesgue$)$ measurable subset. Then,
for each $u\in\mathbb{R,}$
\[
A_u(X,T) = \{t\in T\dvt  X(t)\geq u\}
\]
is called the \textit{excursion set} of $X$ in $T$ over the level $u$.
\end{definition}

Let $\nu_d(B)$ be the volume (i.e., the Lebesgue measure) of a
measurable set $B\subset\mathbb{R}^d$ and~$\ind\{C\}$ denote the
indicator of a set $C$.

Since $X$ is measurable, the volume of the excursion set\vspace*{-1pt}
\[
\nu_d(A_u(X,T)) = \int_{T} \ind\{X(t)\geq u\}\,\mathrm{d}t\vspace*{-1pt}
\]
is a random variable for each $u\in\mathbb{R}$ and any measurable
set $T\subset\mathbb{R}^d$.\vspace*{-2pt}

\vspace*{-3pt}\subsection{Growing sets}\vspace*{-3pt}

Denote the boundary of a set $B\subset\mathbb{R}^d$ by $\partial
B$. The \textit{Minkowski sum} of two sets, $A$, $B\subset\mathbb{R}^d$,
is given by $A\oplus B=\{x+y\dvt x\in A, y\in B\}$. The following
concept of ``regular growth'' for a family of subsets in
$\mathbb{R}^d$ will be used in the sequel.\vspace*{-2pt}

\begin{definition}
A sequence, $(W_n)_{n\in\mathbb{N}}$, of bounded measurable
sets, $W_n\subset\mathbb{R}^d$, tends to infinity in the \textit{Van
Hove sense}
$($VH-growing$)$
if, for any $\epsilon>0$, one has\vspace*{-1pt}
%
\begin{equation}\label{vh}
\nu_d(W_n)\to\infty  \quad \mbox{and} \quad  \frac{\nu_d(\partial W_n
\oplus B_{\epsilon}(0))}{\nu_d(W_n)} \to0\vspace*{-1pt}
\end{equation}
as $n\to\infty$, where $B_{\epsilon}(0)=\{x\in\mathbb{R}^d\dvt  \|
x \|_2\leq
\epsilon\}$ is the closed ball in $\mathbb{R}^d$ with center at the
origin $0\in\mathbb{R}^d$ and radius $\epsilon$.\vspace*{-1pt}
\end{definition}

If $W_n=(a(n),b(n)]=(a_1(n),b_1(n)]\times\cdots \times
(a_d(n),b_d(n)]$ is a parallelepiped, then $W_n \to\infty$ in the
Van Hove sense if and only if $b_k(n)-a_k(n)\to\infty$ as $n\to
\infty$ for $k=1,\ldots,d$.\vspace*{-2pt}

\begin{definition}
A sequence of finite sets $U_n\subset\mathbb{Z}^d$ tends to infinity
in a \textit{regular way} if\vspace*{-1pt}
%
\begin{equation}\label{reg}
\frac{\operatorname{\mathsf{card}} \delta U_n}{\operatorname{\mathsf{card}} U_n}\to0   \qquad \mbox{as }
 n\to\infty;\vspace*{-1pt}
\end{equation}
cf. \eqref{vh}. Here $\delta U_n =\{j \in\mathbb{Z}^d
\setminus U_n\dvt  \dist(j,U_n)=1\}$ and $\dist(j,U_n)= \min_{k\in
U_n}\|j- k\|_{\infty}$.\vspace*{-1pt}
\end{definition}

\section{Central limit theorem}\vspace*{-3pt}\label{sectCLT}

Now we state and prove a CLT for the volume of excursion sets of
random fields. Ivanov and Leonenko \cite{Ivanov} studied
stationary and isotropic Gaussian random fields. In our approach,
the isotropy of Gaussian fields is not required. Moreover, we
consider a more general class of random fields possessing the
quasi-association property. To avoid long formulations, we
introduce the following two conditions for a random field $
X=\{X(t), t\in\mathbb{R}^d\}$.
\begin{longlist}[(B)]
\item[({A})] $X$ is
quasi-associated and strictly stationary such that $X(0)$ has a
bounded density. Assume that the covariance function of $X$ is
continuous and there exists some $\alpha>3d $ such that\vspace*{-1pt}
%
\begin{equation}\label{eq:covariance1}
|\operatorname{\mathsf{cov}}(X(0),X(t))| =
\mathrm{O} (\|t\|_2^{-\alpha} )\qquad
\mbox{as }  \|t\|_2\to\infty .\vspace*{-1pt}
\end{equation}
\item[({B})] $X$ is Gaussian and stationary. Suppose that its
continuous covariance function satisfies \eqref{eq:covariance1} for
some $\alpha>d$.
\end{longlist}

Notice that continuity of the covariance
function of
$X$ implies the existence of a measurable modification of this
field. We consider only such versions of $X$. We exclude the
trivial case when $X(t)=$ const a.s. for all $t\in\mathbb{R}^d$.

\vspace*{-2pt}\subsection{Quasi-associated random fields}\vspace*{-2pt}

To prove the CLT for the volume of excursion sets of a random
field satisfying condition~({A}), we need the following
auxiliary result.\vspace*{-1pt}
\begin{lemma}[(\cite{Bulinski}, Lemma 7.3.4)]\label{lemma:kovarianz1}
Let $\{U,V\}\in\mathsf{QA,}$ where random variables $U$ and $V$ are
square-integrable and have densities bounded by $a>0$. Then\vspace*{-1pt}
\[
| \operatorname{\mathsf{cov}}(\ind\{U\geq u\},\ind\{V\geq v\})|
\leq3\cdot2^{2/3} a^{2/3}| \operatorname{\mathsf{cov}}(U,V)|^{1/3}\vspace*{-1pt}
\]
for arbitrary $u,v\in\mathbb{R}$.\vspace*{-1pt}
\end{lemma}

\begin{theorem}\label{Hauptsatz_general}
Let $ X=\{X(t), t\in\mathbb{R}^d\}$ be a random field satisfying
condition (\textup{{A}}). Then, for any sequence of VH-growing sets
$W_n\subset\mathbb{R}^d $ and each $u\in\mathbb{R}$, one has\vspace*{-1pt}
%
\begin{equation}\label{clt}
\frac{\nu_d(A_{u}(X,W_n)) - \nu_d(W_n) {\mathsf
{P}}(X(0)\geq u) }{\sqrt{\nu_d(W_n)}}
\xrightarrow{d} Y_u\sim\mathcal{N}(0,\sigma^2(u)), \qquad   n\to\infty.\vspace*{-1pt}
\end{equation}
Here $\xrightarrow{d}$ denotes convergence in distribution, $Y_u$ being
a Gaussian random variable with mean zero and variance\vspace*{-1pt}
%
\begin{equation}\label{asdis}
\sigma^2(u)=\int_{\mathbb{R}^d}
\operatorname{\mathsf{cov}}\bigl(\ind\{X(0)\geq u\},\ind\{X(t)\geq u\}\bigr)\,\mathrm{d}t\in\mathbb{R}_+.\vspace*{-1pt}
\end{equation}
\end{theorem}
\begin{pf} Fix any $u \in \mathbb{R}$ and transform a random
field $X$ into a field $Z = \{Z(j),j \in \mathbb{Z}^d\},$ setting\vspace*{-1pt}
%
\begin{equation}\label{def:transform}
Z(j) = \int_{Q_j} \ind\{X(t)\geq u\}\,\mathrm{d}t - {\mathsf
{P}}\bigl(X(0)\geq u\bigr),\vspace*{-1pt}\vadjust{\goodbreak}
\end{equation}
where the unit cubes
\[
Q_j=\{x=(x_1,\ldots,x_d)^\top\in\mathbb{R}^d\dvt 0< x_k\leq1,
 k=1,\ldots ,d\} \oplus\{j\}, \qquad
j=(j_1,\ldots,j_d)^\top\in\mathbb{Z}^d.
\]
The Fubini theorem implies $\mathsf{E}Z(j)=0$ for any $j\in\mathbb{Z}^d$.
It is easily seen that the field $Z$ is strictly stationary and
square-integrable. Introduce
%
\begin{equation}\label{jn}
J_n^-=\{j\in\mathbb{Z}^d\dvt Q_j \subset W_n\}, \qquad
J_n^+=\{j\in\mathbb{Z}^d\dvt Q_j \cap W_n\neq\varnothing\}
\end{equation}
and
\[
W_n^-=\bigcup_{j\in J_n^-} Q_j,  \qquad   W_n^+=\bigcup_{j\in J_n^+} Q_j.
\]
Due to \eqref{vh}, we conclude (see \cite{Bulinski}, Lemma 3.1.2)
that
%
\begin{equation}\label{growth}
\nu_d(W_n^-)\to
\infty   \quad \mbox{and} \quad   \nu_d(W_n^-)/\nu_d(W_n^+)\to
1   \qquad \mbox{as }  n\to\infty.
\end{equation}
Write
%
\begin{eqnarray}\label{eq:sum1}
&&\frac{\nu_d(A_u(X,W_n)) - \nu_d(W_n){\mathsf
{P}}(X(0)\geq
u)}{\sqrt{\nu_d(W_n)}}\nonumber\\[-1pt]
&& \quad =\frac{\nu_d(A_u(X,W_n^-)) - \nu_d(W_n^-){\mathsf
{P}}(X(0)\geq
u)}{\sqrt{\nu_d(W_n)}}\\[-1pt]
&& \qquad {} + \frac{\nu_d(A_u(X,W_n))
-\nu_d(A_u(X,W_n^-))-( \nu_d(W_n)- \nu_d(W_n^-)){\mathsf
{P}}(X(0)\geq u)
}{\sqrt{\nu_d(W_n)}} .
\nonumber
\end{eqnarray}
We prove that the second term on the right-hand side in \eqref
{eq:sum1} tends to zero in
probability. By Chebyshev's inequality, it suffices to show that
\[
\operatorname{\mathsf{var}}\bigl(\nu_d(A_u(X,W_n))- \nu_d(A_u(X,W_n^-))\bigr)/\nu_d(W_n) \to
0, \qquad    n\to\infty.
\]
Set
\[
Y_n(j) =\int_{Q_j\cap W_n} \ind\{X(t)\geq u\}\,\mathrm{d}t - \nu_d(Q_j\cap
W_n) {\mathsf{P}}\bigl(X(0)\geq u\bigr) \qquad\mbox{for } j\in
\mathbb{Z}^d,\   n\in\mathbb{N} .
\]
Note that $Y_n(j)=Z(j)$ for $j\in J_n^-$ and $n\in\mathbb{N}$
(clearly $Y_n(j)$ and $Z(j)$ depend on $u$ as well). Applying the
Fubini theorem and Lemma \ref{lemma:kovarianz1}, we get
%
\begin{eqnarray} \label{dgr}
&&\operatorname{\mathsf{var}}\bigl(\nu_d(A_u  (X,W_n))-
\nu_d(A_u(X,W_n^-))\bigr)\nonumber\\[-2pt]
 && \quad =\operatorname{\mathsf{var}} \biggl(\sum_{j\in J_n^+\setminus J_n^-} Y_n(j) \biggr)
\nonumber\\[-2pt]
&& \quad  \leq\sum_{j,m\in J_n^+\setminus J_n^-} \int_{Q_j\times
Q_m}\bigl|\operatorname{\mathsf{cov}}\bigl(\ind\{X(s)\geq u\},\ind\{X(t)\geq
u\}\bigr)\bigr|\,\mathrm{d}s\,\mathrm{d}t
\\[-2pt]
&& \quad  \leq\nu_d(W_n^+ \setminus W_n^-) \sum_{j\in\mathbb{Z}^d}
\int_{Q_0\times Q_j} C_1|\operatorname{\mathsf{cov}}(X(s),X(t))|^{1/3}\,\mathrm{d}s\,\mathrm{d}t\nonumber\\[-2pt]
&& \quad  \leq\nu_d(W_n^+ \setminus W_n^-) \Biggl(C_2 + C_3
\sum_{r=r_0}^{\infty}
\sum_{j\in\mathbb{Z}^d\dvt\|j\|_{\infty}=r}
\int_{Q_0\times Q_j} \|s-t\|_2^{-\alpha/3}\,\mathrm{d}s\,\mathrm{d}t \Biggr)\nonumber\\[-2pt]
&& \quad
\leq\nu_d(W_n^+ \setminus W_n^-) \Biggl(C_2 + C_4\sum_{r=1}^{\infty}
r^{d-1} r^{-\alpha/3} \Biggr) = C_5\nu_d(W_n^+ \setminus W_n^-)
\nonumber
\end{eqnarray}
for some $r_0>0$ and all $n\in\mathbb{N}$. The factors $C_i$ do not
depend on $n$. We used the
inequality $|\operatorname{\mathsf{cov}}(X(s),X(t))|\leq\tau^2$ for all
$s,t\in\mathbb{R}^d $, which is satisfied as $\operatorname{\mathsf{var}}X(t)=\tau^2$
for any $t\in\mathbb{R}^d$. We also took into account that
\[
\operatorname{\mathsf{card}}\{j\in\mathbb{Z}^d\dvt\|j\|_{\infty}=r\} \leq C_6 r^{d-1}
\]
for each $r\in\mathbb{N}$ and employed the inequality
\eqref{eq:covariance1} with $\alpha>3d$.

By \eqref{growth}, \eqref{dgr} and in view of the relation $\nu
_d(W_n^-) \leq\nu_d(W_n) \leq\nu_d(W_n^+),$
we get
\[
\operatorname{\mathsf{var}} \bigl(\bigl(\nu_d(A_u(X,W_n))-
\nu_d(A_u(X,W_n^-))\bigr)/\sqrt{\nu_d(W_n)} \bigr) \to0, \qquad    n\to
\infty.
\]
Now we show that $J_n^-$, introduced in \eqref{jn}, tends to
infinity in a regular way as $n\to\infty$. Indeed, $J_n^- \subset
J_n \subset J_n^+$, where $J_n:= W_n \cap\mathbb{Z}^d$,
$n\in\mathbb{Z}^d$. Due to \cite{Bulinski}, Lemma 3.1.5, $J_n$
tends to infinity in a regular way. Thus, it suffices to mention
that $\delta J_n^- \subset\delta J_n \cup(J_n \setminus J_n^-)$
and apply the relations $\operatorname{\mathsf{card}} \delta J_n /\operatorname{\mathsf{card}}  J_n
\to0$ and $\operatorname{\mathsf{card}}  J_n^+ /\operatorname{\mathsf{card}} J_n^- \to1$ as $n\to
\infty$. Lemma~3.1.6 of~\cite{Bulinski} implies that $W_n^- =
\bigcup_{j\in J_n^-}Q_j$ tends to infinity in the Van Hove sense as
$n\to\infty$.

So, while establishing \eqref{clt}, we can assume w.l.g. that
$W_n=W_n^-$, that is, $W_n$ is a finite union of cubes $Q_j$
$(n\in\mathbb{N})$ and the sequence $(W_n)_{n\in\mathbb{N}}$ is
VH-growing.

Observe that
\[
\frac{\nu_d(A_u(X,W_n)) - \nu_d(W_n){\mathsf
{P}}(X(0)\geq u)}{\sqrt{\nu_d(W_n)}}
= \frac{\sum_{j\in W_n \cap\mathbb{Z}^d} Z(j)}{\sqrt{\nu_d(W_n)}}:=S_n.
\]
As $ X=\{X(t), t\in\mathbb{R}^d\}\in\mathsf{QA}$, it follows from
\eqref{int} and \eqref{eq:covariance1} that $X \in(\mathit{BL},\theta)$
with
\[
\theta_r(X) =\mathrm{O}(r^{-\alpha+d})    \qquad \mbox{as }   r\to\infty  (r>0).
\]
For $\gamma>0$ (and $u$ fixed) introduce the Lipschitz functions
$h_{\gamma}\dvtx\mathbb{R}\rightarrow\mathbb{R}$ by the
formula
%
\begin{equation}\label{fh}
h_{\gamma}(x) =
\cases{\displaystyle
0 ,&\quad  if $x \leq u-\gamma$, \cr\displaystyle
(x-u+\gamma)/\gamma,&\quad  if $u-\gamma< x \leq u$,\cr\displaystyle
1 ,&\quad  otherwise.
}
\end{equation}
Superposition of two Lipschitz functions is also a Lipschitz one.
Thus, for $n\in\mathbb{N}$ and $\gamma>0$, the random field
$Z_{n,\gamma} = \{Z_{n,\gamma}(j), j=(j_1,\ldots,j_d)^\top\in
\mathbb{Z}^d\}\in
(\mathit{BL},\theta),$ where
%
\begin{equation}\label{zn}
Z_{n,\gamma}(j) = \frac{1}{n^d}\sum_{k_1,\ldots,k_d = 1}^{n}
h_{\gamma} \biggl(X \biggl(j_1 + \frac{k_1}{n},\ldots,j_d +
\frac{k_d}{n} \biggr) \biggr)-\E h_\gamma(X(0))\vadjust{\goodbreak}
\end{equation}
and the terms of a sequence $\theta(Z_{n,\gamma})$ admit the estimate
%
\begin{equation}\label{bldep}
\theta_r(Z_{n,\gamma})\leq C_7 \gamma^{-2}
r^{-\alpha+d}, \qquad    r\in\mathbb{N},
\end{equation}
with $C_7$ depending neither on $\gamma$ nor on $n$.

It is not difficult to verify that the finite-dimensional distributions
of the fields $Z_{n,\gamma}$ weakly converge to the corresponding
ones of the field $Z_{\gamma}$ as $n\to\infty,$ where
%
\begin{equation}\label{z}
Z_{\gamma}(j) = \int_{Q_j} h_{\gamma}(X(t))\,\mathrm{d}t - \E h_\gamma
(X(0)),  \qquad  j\in\mathbb{Z}^d.
\end{equation}
Consequently (see \cite{Bulinski}, Lemma 1.5.16), we can claim that
$Z_{\gamma}\in(\mathit{BL},\theta)$ and $\theta_r(Z_{\gamma})$ is bounded
by the right-hand side of inequality \eqref{bldep}. Theorem
3.1.12 of \cite{Bulinski}, guarantees that, for each $\gamma
>0$,
%
\begin{equation}\label{auxclt1}
S_n(\gamma):=\frac{\sum_{j\in W_n\cap\mathbb{Z}^d}Z_{\gamma}(j)}
{\sqrt{\nu_d(W_n)}} \stackrel{d}\to Y_{u,\gamma}\sim
\mathcal{N}(0,\sigma^2(u,\gamma)),  \qquad   n\to\infty,
\end{equation}
where
\[
\sigma^2(u,\gamma)= \sum_{j\in\mathbb{Z}^d}\operatorname{\mathsf{cov}}(Z_{\gamma}(0),
Z_{\gamma}(j)) = \int_{\mathbb{R}^d}\operatorname{\mathsf{cov}}(h_{\gamma}(X(0)),
h_{\gamma}(X(t)))\,\mathrm{d}t \in\mathbb{R}_+.
\]
Therefore, to prove \eqref{clt}, two steps remain. First
of all, we estimate the difference of the characteristic functions
of the random variables $S_n(\gamma)$ and $S_n$ and show that it
tends to zero as $\gamma\to0+$. After that, we verify that
%
\begin{equation}\label{est3}
\sigma^2(u,\gamma)\to\sigma^2(u)   \qquad \mbox{as }  \gamma\to0+.
\end{equation}
Set $h(x)=\mathbb{I}\{x\geq u\}$ and
$H_{\gamma}(x)=h_{\gamma}(x)-h(x),$ where $x\in\mathbb{R}$ (and
$u\in\mathbb{R}$ is fixed). Then, for each $\lambda\in\mathbb{R}$,
one has
%
\begin{equation}\label{chf}
\bigl|\mathsf{E}\mathrm{e}^{\mathrm{i}\lambda S_n(\gamma)} - \mathsf{E}\mathrm{e}^{\mathrm{i}\lambda S_n}\bigr| \leq
|\lambda|  \mathsf{E}|S_n(\gamma)-S_n|\leq
|\lambda| \biggl(\frac{V_n(\gamma)}{\nu_d(W_n)} \biggr)^{1/2},
\end{equation}
where $\mathrm{i}^2=-1$ and
\[
V_n(\gamma)=\mathsf{E} \biggl(\sum_{j\in W_n\cap\mathbb{Z}^d}
\int_{Q_j}\bigl(H_{\gamma}(X(t))- \mathsf{E}H_{\gamma}(X(t))\bigr)\,\mathrm{d}t \biggr)^{  2} .
\]
It is easily seen that
%
\begin{equation}\label{est1}
V_n(\gamma) \leq\nu_d(W_n) \int_{\mathbb{R}^d}|\operatorname{\mathsf{cov}}(H_{\gamma}(X(0)),H_{\gamma}(X(t)))|\,\mathrm{d}t.
\end{equation}
Furthermore, we have
\[
|\operatorname{\mathsf{cov}}(H_{\gamma}(X(0)),H_{\gamma}(X(t)))| \leq (\mathsf{E}(H_{\gamma}(X(0)))^2\mathsf{E}(H_{\gamma}(X(t)))^2 )^{1/2}
\leq a \gamma,\vadjust{\goodbreak}
\]
where $a$ is a constant that bounds the density of $X(0)$. If $
|\operatorname{\mathsf{cov}}(X(0),X(t))|^{1/3}\leq\gamma$, then reasoning
similar to that proving Lemma \ref{lemma:kovarianz1} leads to
the inequality
%
\begin{equation}\label{ec}
|\operatorname{\mathsf{cov}}(H_{\gamma}(X(0)),H_{\gamma}(X(t)))|
\leq C(a) |\operatorname{\mathsf{cov}}(X(0),X(t))|^{1/3}
\end{equation}
with some $C(a)>0$. Write $\alpha= 3(d+\mu)$, $\mu
>0$, and take
$R= c\gamma^{- {1}/({d+\mu})}$, where $c>0$. Then,
in view of \eqref{eq:covariance1} and due to the appropriate choice
of $c$, one can conclude that for all $\gamma>0$ small enough,
%
\begin{eqnarray} \label{est2}
F(\gamma)&:=&\int_{\mathbb{R}^d}|\operatorname{\mathsf{cov}}(H_{\gamma}(X(0)),H_{\gamma}(X(t)))|\,\mathrm{d}t
\nonumber
\\[-8pt]
\\[-8pt]
&\leq& a\gamma\omega_d R^d + C_8 \int_{\|t\|_{2} \geq R}
\|t\|_{2}^{-\alpha/3}\,\mathrm{d}t \leq C_9\gamma^{ {\mu}/({d+\mu})},
\nonumber
\end{eqnarray}
where $\omega_d=\pi^{d/2}/ (\Gamma(d/2+1) )$ is the volume
of the unit ball in $\mathbb{R}^d$ with the Euclidean norm.
Consequently, inequalities \eqref{chf}, \eqref{est1} and \eqref
{est2} imply
that the laws of $S_n(\gamma)$ and~$S_n$ are close for all $n$
large enough if $\gamma>0$ is small enough.

Next, we proceed to \eqref{est3}. By the arguments leading
to \eqref{est2} and invoking Lemma \ref{lemma:kovarianz1}, we
deduce that $\sigma^2(u) <\infty$. Similar to \eqref{ec}, one
shows that if $ |\operatorname{\mathsf{cov}}(X(s),X(t))|^{1/3}\leq
\gamma$,~then
%
\begin{equation}\label{est5}
|\operatorname{\mathsf{cov}}(h(X(s)),H_{\gamma}(X(t)))|
\leq D(a) |\operatorname{\mathsf{cov}}(X(s),X(t))|^{1/3}
\end{equation}
with $D(a)>0$ depending on $a$ only. The absolute value of $\sigma
^2(u,\gamma) - \sigma^2(u)$ does not
exceed the following expression:
\[\label{est5a}
F(\gamma) + \int_{\mathbb{R}^d}|\operatorname{\mathsf{cov}}(h(X(0)),H_{\gamma}(X(t)))|\,\mathrm{d}t +
\int_{\mathbb{R}^d}|\operatorname{\mathsf{cov}}(H_{\gamma}(X(0)),h(X(t)))|\,\mathrm{d}t.
\]
Taking into account the above upper bound and relations
\eqref{est2} and \eqref{est5}, we complete the proof of
\eqref{est3}. The asymptotic (finite) variances
$\sigma^2(u,\gamma)$ are non-negative, whence one concludes
that $\sigma^2(u)\geq0$.

In view of \eqref{auxclt1}--\eqref{chf}, the proof
is complete.
\end{pf}

Now we turn to the multidimensional CLT for random vectors,
%
\begin{equation}\label{eq:partialsum}
S_{\vec{u}}(X,W_n) =  (\nu_d(A_{u_1}(X,W_n)),\ldots,\nu
_d(A_{u_r}(X,W_n)) )^{\top}, \qquad    n\in\mathbb{N},
\end{equation}
where $\vec{u}=(u_1,\ldots,u_r)^\top\in\mathbb{R}^r$.

%
\begin{theorem}\label{th:multiclt1}
Let $X=\{X(t), t\in\mathbb{R}^d\}$ be a random field satisfying condition
{\textup{(A)}}. Then, for each $\vec{u}=(u_1,\ldots,u_r)^\top\in\mathbb
{R}^r$ and
any VH-growing sequence $(W_n)_{n\in\mathbb{N}}$ of
subsets of $\mathbb{R}^d$, one has
%
\begin{equation}\label{multclt}
\nu_d(W_n)^{-1/2}\bigl(S_{\vec{u}}(X,W_n)-\nu_d(W_n) P(\vec{u})\bigr)
\xrightarrow{d} V_{\vec{u}}\sim\mathcal{N}(0,\Sigma(\vec{u}))
 \qquad \mbox{as }  n\to\infty,
\end{equation}
where
\[
P(\vec{u}) = \bigl({\mathsf{P}}\bigl(X(0)\geq u_1\bigr),\ldots
,{\mathsf{P}}\bigl(X(0)\geq
u_r\bigr)\bigr)^{\top}\vadjust{\goodbreak}
\]
and $\Sigma(\vec{u})=(\sigma_{lm}(\vec{u}))_{l,m=1}^r$ is an
$(r\times r)$-matrix having the elements
%
\begin{equation}\label{eq:sigmalm}
\sigma_{lm}(\vec{u}) =
\int_{\mathbb{R}^d} \operatorname{\mathsf{cov}}\bigl(\ind\{X(0)\geq u_l\},\ind\{X(t)\geq
u_m\}\bigr)\,\mathrm{d}t.
\end{equation}

\end{theorem}

\begin{pf}
Observe that the convergence of all $r^2$ integrals in (\ref
{eq:sigmalm}) is proved in the same way as that of the integral
representing $\sigma^2(u)$ in the one-dimensional case. The result
follows by using the Cram\'{e}r--Wold device. We omit further details
that are quite similar to those in the proof of Theorem~\ref
{Hauptsatz_general}.
\end{pf}

The last theorem entails:
\begin{corollary}
Let $X=\{X(t), t\in\mathbb{R}^d\}$ be a random field satisfying condition
{\textup{(A)}}. Assume that
$\Sigma(\vec{u})$ is non-degenerate for some
$\vec{u}\in\mathbb{R}^r$. Then, for this $\vec{u}$
and any sequence $(W_n)_{n\in\mathbb{N}}$ of VH-growing
subsets of $\mathbb{R}^d$, one has

%
\[
\nu_d(W_n)^{-1/2}\Sigma(\vec{u})^{-1/2}\bigl(S_{\vec{u}}(X,W_n)-\nu
_d(W_n) P(\vec{u})\bigr)
\xrightarrow{d} V\sim\mathcal{N}(0,\mathrm{I})   \qquad \mbox{as }
n\to\infty;
\]
here {I} denotes the unit $(r\times r)$-matrix.
\end{corollary}

\subsection{Shot noise random fields} \label{sect:shotnoise}

We verify the conditions of Theorem \ref{Hauptsatz_general} for
shot noise random fields. These fields appear naturally in the
theory of disordered structures. Let $\mathcal{B}(\R^d)$
(resp., $\mathcal{B}_0(\R^d)$) be the family of all (bounded)
Borel sets in $\R^d$. A shot noise random field $X=\{X(t), t\in\R^d
\}$ is defined by the relation
\[
X(t)=\sum_{i\in\N} \xi_i \varphi(t-x_i),
\]
where $\{\xi_i\}$ is a family of i.i.d.
non-negative random variables and $\{x_i\}$ is a
homogeneous Poisson point process in $\R^d$ with intensity $\lambda
\in(0,\infty)$, that is, $\{x_i\}$ is the support set of a~random
Poisson counting measure $\{N_B, B\in\mathcal{B}(\R^d)\},$ where
$N_B=\#\{i\dvt x_i\in B\}$ has the following properties:
%
\begin{longlist}[(ii)]
\item[(i)]$N_{B_1},N_{B_2},\ldots$ are independent for pairwise disjoint
$B_1,B_2,\ldots\in\mathcal{B}_0(\R^d)$,
\item[(ii)]$N_B\sim \operatorname{Pois}(\lambda\nu_d(B))$ for all
$B\in\mathcal{B}_0(\R^d)$.
\end{longlist}
Suppose that $\{\xi_i\}$, $N_{(\cdot)}$ are independent, $\E\xi
_i^2<\infty$ and
$\varphi\dvtx\R^d\to\R_+$ is a Borel function.\vadjust{\goodbreak}

For the shot-noise field $X$ introduced above, we impose the
condition:

({C}) $X(0)$ has a bounded density and
for a function $\varphi$ bounded and uniformly continuous on $\R^d,$
%
\begin{equation}\label{eq:shotnoise}
\varphi(t)\leq\varphi_0(\|t\|_2) = \mathrm{O} (\|t\|_2^{-\alpha} )
 \qquad
\mbox{as }  \|t\|_2\to\infty,
\end{equation}
where $\alpha>3d$ and $\varphi_0\dvtx\R_+\to\R_+$.

\begin{proposition}
The statement of Theorem $\ref{Hauptsatz_general}$ holds for a
random field $X$ satisfying condition (\textup{{C}}).
\end{proposition}
\begin{pf}
By \cite{Bulinski}, Theorem 1.3.8, $X$ is associated and hence
quasi-associated. Moreover, it is strictly stationary with covariance
function given, for example, in \cite{Bulinski}, Theorem 2.3.6.
The continuity of the covariance function follows from
the inequality
\[
 |\cov(X(0),X(s)) - \cov(X(0),X(t)) | \leq\lambda \mathsf{E} \xi_1^2 \sup_{y\in\R^d} |\varphi(t-y)-\varphi(s-y)| \int_{\R^d}
|\varphi(y)|\,\mathrm{d}y
\]
and the uniform continuity of $\varphi$. Corollary 2.3.7 of
\cite{Bulinski} yields the desired bound for the covariance
function in condition ({A}). The proof is complete.
\end{pf}

Note that
the characteristic function of
$X(0)$, provided by \cite{Bulinski}, Lemma 1.3.7, is integrable if
%
\begin{equation}\label{eq:BoundedDens}
\int_{\R}  \biggl|  \exp \biggl\{\lambda\int_{\R^d}
 \bigl(\varphi_\xi(s\varphi(t))-1 \bigr)\,\mathrm{d}t \biggr\}  \biggr|\,\mathrm{d}s<\infty.
\end{equation}
Thus, \eqref{eq:BoundedDens} guarantees the
existence of the bounded density of $X(0)$.



Condition \eqref{eq:BoundedDens} can be easily verified in a
number of special cases; for instance, if $\xi_1=\operatorname{const}>0$
a.s. and $\varphi(t)=a\exp\{-b \|t\|_2 \}$ or $\varphi(t)=a\min\{
1, \|t\|_2^{-b} \}$ with $a,b>0$.

\subsection{Gaussian random fields} \label{sect:Gaussian}

In contrast to Lemma \ref{lemma:kovarianz1}, we obtain a sharper
estimate for the covariance of indicator functions in the Gaussian
case. Our result extends formula (2.7.1) of \cite{Ivanov}.
Let $\Phi$ and~$\Psi$ stand for the cumulative distribution
function and the tail distribution function of a~standard Gaussian
random variable, respectively.

\begin{lemma}\label{lemma:kovarianz2}
Let $ (U, V)^\top$ be a Gaussian random vector in $\mathbb{R}^2$
such that $ U \sim\mathcal{N}(a,\tau^2) $, $ V \sim
\mathcal{N}(a,\tau^2), $ where $a\in\mathbb{R}$, $\tau>0$ and
correlation coefficient $
\corr (U,V ) = \rho$. Then, for any $u,v\in\R$ and
$\rho\in(-1,1)$, the following equality holds:
%
\begin{eqnarray}\label{eq:GaussCov}
&&\cov(\ind\{U\geq u\},\ind\{V\geq v\})\nonumber
\\[-8pt]
\\[-8pt]
&& \quad =\frac{1}{2\uppi}\int_0^{\rho}
\frac{1}{\sqrt{1-r^2}} \exp \biggl\{-\frac
{(u-a)^2-2r(u-a)(v-a)+(v-a)^2}{2\tau^2(1-r^2)} \biggr\}\,\mathrm{d}r .
\nonumber
\end{eqnarray}
In particular, for $u=v$, one has
\[
\cov(\ind\{U\geq u\},\ind\{V\geq u\})=
\frac{1}{2\uppi} \int_0^{\rho}
\frac{1}{\sqrt{1-r^2}}\exp \biggl\{-\frac{(u-a)^2}{\tau
^2(1+r)} \biggr\}\,\mathrm{d}r.
\]
Moreover, for any $u,v\in\R$ and $\rho\in[-1,1]$, one has the inequality
%
\begin{equation}\label{eq:boundGauss}
 |\cov(\ind\{U\geq u\},\ind\{V\geq v\}) | \leq|\rho|/4.
\end{equation}
\end{lemma}
%
%
\begin{pf}
Using the transformation $x \mapsto(x-a)/\tau$, $x\in\mathbb{R,}$
we can assume w.l.g. that $U
\sim\mathcal{N}(0,1) $ and $ V \sim\mathcal{N}(0,1)$.
Let $\rho\in(-1,1)$. The probability density
\[
f_{U,V}(x,y)= \frac{1}{2\uppi\sqrt{1-\rho^2}}
\exp \biggl\{-\frac{x^2-2\rho xy+y^2}{2(1-\rho^2)} \biggr\}
\]
of the bivariate Gaussian random variable $(U,V)^\top$ is
invariant under the transformation $x\mapsto-x$ and $y\mapsto
-y$, $(x,y)^\top\in\R^2$. Therefore,
\[
\cov(\ind\{U\geq u\},\ind\{V\geq v\})=\cov(\ind\{U\leq-u\},\ind\{
V\leq-v\}), \qquad
u,v\in\R.
\]
It is well known (see, e.g., \cite{Cramer}, formulae (21.12.5) and
(21.12.6)) that
\[
f_{U,V}(x,y) = \sum_{k=0}^{\infty}
\frac{\Phi^{ (k+1 )}(x)
\Phi^{ (k+1 )}(y)}{k!}\rho^k, \qquad    x,y\in\mathbb{R},
\]
where $\Phi^{(k)}(x)=\mathrm{d}^k\Phi(x)/\mathrm{d}x^k$
and, for any $u,v\in\mathbb{R}$,
\[
\int_{-\infty}^u\int_{-\infty}^v f_{U,V}(x,y)\,\mathrm{d}x\,\mathrm{d}y = \sum
_{k=0}^{\infty}
\frac{\Phi^{(k)}(u)\Phi^{(k)}(v)}{k!}\rho^k.
\]
Hence, for each $u,v\in\mathbb{R}$,
\begin{eqnarray*}
\mathsf{E}  \ind\{U\leq-u\}\ind\{V\leq-v\}& =&
\int_{-\infty}^{-u} \int_{-\infty}^{-v} f_{U,V}(x,y)\,\mathrm{d}x\,\mathrm{d}y
= \sum_{k=0}^{\infty}
\frac{\Phi^{(k)}(-u)\Phi^{(k)}(-v)}{k!}\rho^k\\
&=&
\Phi(-u)\Phi(-v)+\sum_{k=1}^{\infty}
\frac{\Phi^{(k)}(-u)\Phi^{(k)}(-v)}{k!}\rho^{k}\\
&=& \Phi(-u)\Phi(-v)+\int_0^{\rho}
\sum_{k=0}^{\infty}\frac{\Phi^{(k+1)}(-u)\Phi^{(k+1)}(-v)}{k!}r^k\,\mathrm{d}r\\
&=& \Phi(-u)\Phi(-v)+\int_0^{\rho}f_{U(r),V(r)}(u,v)\,\mathrm{d}r,
\end{eqnarray*}
where $(U(r),V(r))^\top$ is a centered bivariate Gaussian vector with
$\mathsf{E}U(r)^2 = \mathsf{E}V(r)^2 =  1$ and $\operatorname{\mathsf{cov}}(U(r),V(r))=r$.
Consequently, we get
\begin{eqnarray*}
\operatorname{\mathsf{cov}}(\ind\{U\leq-u\},\ind\{V\leq-v\}) &=&
\int_0^{\rho}f_{U(r),V(r)}(u,v)\,\mathrm{d}r\\
&=& \frac{1}{2\uppi} \int_0^{\rho}
\frac{1}{\sqrt{1-r^2}} \mathrm{e}^{-({u^2-2ruv+v^2})/({2(1-r^2)})}\,\mathrm{d}r.
\end{eqnarray*}
Passing to random variables $U$ and $V$ with arbitrary mean $a$
and variance $\tau^2>0$ gives the formula \eqref{eq:GaussCov}.

To prove inequality \eqref{eq:boundGauss} for $\rho\in(-1,1)$,
write
\[
 |\operatorname{\mathsf{cov}}(\ind\{U\leq-u\},\ind\{V\leq-v\}) | \leq
\frac{1}{2\uppi} \int_0^{|\rho|}
\frac{1}{\sqrt{1-r^2}}\,\mathrm{d}r\leq\frac{1}{2\uppi} \arcsin|\rho|
\]
and notice that $\arcsin|\rho|\leq\uppi  |\rho|/2 $.

The case $|\rho|=1$ is trivial, as $|\mathsf{P}(A\cap B)-\mathsf{P}(A)
\mathsf{P}(B)|\le1/4$ for any $A,B\in\mathcal{F}$.
\end{pf}
The following result generalizes the corresponding one established in
\cite{Ivanov} (see page 80),
where the isotropy of the Gaussian random field was assumed.
A central limit theorem for nonlinear transformations of a homogeneous
Gaussian random field
was used there.\vspace*{-3pt}
%
\begin{theorem}\label{th:gaussian}
Let $ X=\{X(t), t\in\R^d\} $ be a Gaussian stationary random
field satisfying condition {\textup{(B)}} and $X(0)\sim\mathcal{N}(a,\tau
^2)$. Then, for each $u\in\R$ and any sequence of $VH$-growing sets
$W_n\subset\R^d $, one has
\[
\frac{\nu_d(A_u(X,W_n)) - \nu_d(W_n) \Psi((u-a)/\tau)}{\sqrt{\nu
_d(W_n)}} \xrightarrow{d} Y_u\sim
\mathcal{N}(0,\sigma^2(u))
\]
as $ n\to\infty$. The variance $\sigma^2(u)$ introduced in \eqref{asdis}
can be written in the following form:
%
\begin{equation}\label{dispers}
\sigma^2(u) = \frac{1}{2\uppi}\int_{\R^d} \int_0^{\rho(t)}
\frac{1}{\sqrt{1-r^2}}\exp \biggl\{-\frac{(u-a)^2}{\tau
^2(1+r)} \biggr\}\,\mathrm{d}r\,\mathrm{d}t,
\end{equation}
where $\rho(t)=\operatorname{\mathsf{corr}}(X(0),X(t))$. In particular,
\[
\sigma^2(a) = \frac{1}{2\uppi}\int_{\R^d} \arcsin(\rho(t))\,\mathrm{d}t.\vspace*{-3pt}
\]
\end{theorem}
\begin{pf}
For the Gaussian field $X,$ we have ${\mathsf{P}}(X(0)\geq
u)=\Psi((u-a)/\tau)$. Now we apply the upper bound
\eqref{eq:boundGauss} to obtain $|\operatorname{\mathsf{cov}}(X(0),X(t))|$ instead
of $|\operatorname{\mathsf{cov}}(X(0),X(t))|^{1/3}$ in the estimates used in the
proof of Theorem \ref{Hauptsatz_general}. This leads to the
hypothesis that $\alpha
>d$ in~\eqref{eq:covariance1}, whereas in condition ({A}) we
assumed $\alpha>3d$.
Note that Gaussian fields are quasi-associated~\cite{Shash}.

Finally, we express $\sigma^2(u)$ (see \eqref{asdis}) in terms of the
covariance
function of $X$ as in the proof of Lemma \ref{lemma:kovarianz2}, and
this yields \eqref{dispers}. The proof is complete.\vspace*{-3pt}
\end{pf}
%
%
\begin{theorem}\label{th:multiclt2}
Let $X=\{X(t), t\in\R^d\}$ be a random field satisfying condition
{\textup{(B)}} and $X(0)\sim\mathcal{N}(a,\tau^2)$. Then, for each $\vec
{u}=(u_1,\ldots,u_r)^\top\in\R^r$ and any sequence
$(W_n)_{n\in\N}$ of VH-growing subsets of $\mathbb{R}^d$, one has
%
\begin{equation}\label{eq:multiclt2}
\nu_d(W_n)^{-1/2}\bigl(S_{\vec{u}}(X,W_n)-\nu_d(W_n) \Psi(\vec{u})\bigr)
\xrightarrow{d}V_{\vec{u}}\sim\mathcal{N}(0,\Sigma(\vec{u}))
   \qquad \mbox{as }  n\to\infty.
\end{equation}
Here, $ \Psi(\vec{u})=(\Psi((u_1-a)/\tau),\ldots,\Psi
((u_r-a)/\tau))^{\top}$ and
$\Sigma(\vec{u})= (\sigma_{lm}(\vec{u}))_{l,m=1}^r$ is a matrix
having the elements
%
\begin{equation}\label{elm}
\sigma_{lm}(\vec{u}) = \frac{1}{2\uppi}\int_{\R^d} \int_0^{\rho
(t)}g(r)\,\mathrm{d}r\,\mathrm{d}t,\vadjust{\goodbreak}
\end{equation}
where
\[
g(r) = \frac{1}{\sqrt{1-r^2}} \exp \biggl\{-\frac
{(u_l-a)^2-2r(u_l-a)(u_m-a)+(u_m-a)^2}{2\tau^2(1-r^2)} \biggr\}
\]
and $\rho(t) = \corr(X(0),X(t))$. If $\Sigma(\vec{u})$ is
non-degenerate, we obtain by virtue of
$\eqref{eq:multiclt2}$
\[
\nu_d(W_n)^{-1/2}\Sigma(\vec{u})^{-1/2}\bigl(S_{\vec{u}}(X,W_n)-\nu
_d(W_n) \Psi(\vec{u})\bigr)
\xrightarrow{d} \mathcal{N}(0,\mathrm{I}), \qquad   n\to\infty,
\]
where {I} is the unit $(r\times r)$-matrix.
\end{theorem}
\begin{pf} Employing Lemma \ref{lemma:kovarianz2}, one can repeat
the reasoning proving Theorem \ref{th:multiclt1}. Clearly,
${\mathsf{P}}(X(0)\geq u_l) = \Psi((u_l-a)/\tau)$,
$l=1,\ldots,r$. The matrix elements $\sigma_{lm}(\vec{u})$ for
$l,m=1,\ldots,r$ can be calculated by way of
\eqref{eq:GaussCov}.
\end{pf}

\subsubsection*{Formulae in the isotropic case}

For the isotropic case, we use the change of variables (passing from
$t=(t_1,\ldots,t_d)^\top$ to spherical coordinates) in integrals
\eqref{dispers} and \eqref{elm} to obtain the following statement.

\begin{corollary}
Let the random field $X=\{X(t), t\in\R^d\}$ satisfying the
conditions of Theorem $\ref{th:gaussian}$ be isotropic $(d\geq2)$. Then
\[
\sigma^2(u) = \frac{\mathrm{d}\omega_d}
{2\uppi} \int_0^{\infty} v^{d-1}
\int_0^{\rho(v)} \frac{1}{\sqrt{1-r^2}}
\exp\biggl \{-\frac{(u-a)^2}{\tau^2(1+r)} \biggr\}\,\mathrm{d}r\,\mathrm{d}v,
\]
where $\rho(v)=\corr(X(0),X(t))$ if $|t|=v$.
For $u=a,$ one has
\[
\sigma^2(a) = \frac{\mathrm{d}\omega_d}
{2\uppi}
\int_0^{\infty} v^{d-1}\arcsin(\rho(v))\,\mathrm{d}v.
\]
In the multivariate case, \eqref{elm} can be written as follows:
\[
\sigma_{lm}(\vec{u}) = \frac{\mathrm{d}\omega_d}{2\uppi}\int_0^{\infty} v^{d-1}
\int_0^{\rho(v)} g(r)\,\mathrm{d}r\,\mathrm{d}v
\]
for $l,m=1,\ldots,r$.
\end{corollary}


\section{Statistical version of the CLT}\label{sectStatCLT}

Now we provide a statistical version of the CLT involving random
self-normalization. Let $r\in\N$ be the
number of levels to observe.
\begin{theorem}\label{th:statvers}
Let $X=\{X(t), t\in\R^d\}$ be a random field
satisfying condition {\textup{(A)}}.
Let $u_k\in\R$, $k=1,\ldots,r$ and
$(W_n)_{n\in\N}$ be a sequence of VH-growing sets. Furthermore,
let $\hat{C}_n =
 (\hat{c}_{nlm} )_{l,m=1}^r$
be statistical estimates for non-degenerate
asymptotic covariance
matrix~$\Sigma$ with elements $\sigma_{lm}$
given by \eqref{eq:sigmalm}.
Assume that $\hat{c}_{nlm} \xrightarrow{p} \sigma_{lm}$
as $n\to\infty$ for any $l,m=1,\ldots,r,$
where $\xrightarrow{p}$ denotes convergence in probability. Then
\[
\hat{C}_n^{-1/2}\nu_d(W_n)^{-1/2}\bigl(S(W_n)- \nu_d(W_n) P(\vec{u})\bigr)
\xrightarrow{d} \mathcal{N}(0,\mathrm{I})   \qquad \mbox{as }
n\to\infty.\vspace*{-2pt}
\]
\end{theorem}
\begin{pf}
It suffices to use Theorem \ref{th:multiclt1} and elementary
properties of
the convergence in probability and in law for random vectors.\vspace*{-2pt}
\end{pf}

One feasible estimator for the asymptotic covariance matrix $\Sigma$ that
arose in the multivariate CLT, see Theorem \ref{th:multiclt1}, can be
called a \textit{subwindow estimator} \cite{Pantle} and is constructed
as follows.
Let $(V_n)_{n\in\N}$ and $(W_n)_{n\in\N}$ be
sequences of VH-growing sets (not necessarily rectangles) such
that $V_n\subset W_n$, $n\geq1$. Consider $N(n)$ subwindows
$V_{n,1},\ldots,V_{n,N(n)}$, where $(N(n))_{n\in\N}$ is an
increasing sequence of integers with $\lim_{n\to\infty}
N(n)=\infty$, and $V_{n,j}=V_n\oplus\{h_{n,j}\}$ are subwindows
that are translated by certain vectors $h_{n,j}\in\R^d$,
$j=1,\ldots,N(n)$. Assume that $\bigcup_{j=1}^{N(n)} V_{n,j}\subseteq
W_n$ for each $n\in\N$ and there exists some $r>0$ such that\vspace*{-2pt}
\[
V_{n,j}\cap V_{n,i} \subset\partial V_{n,j} \oplus B_r(0)   \qquad \mbox
{for } i,j\in\{1,\ldots,N(n)\}    \mbox{ with } i\neq j.\vspace*{-2pt}
\]
Denote by\vspace*{-2pt}
\[
\hat{\mu}^{(j)}_{nk}=\frac{1}{\nu_d(V_n)} \int_{V_{n,j}} \ind\{
X(t)\geq u_k\}\,\mathrm{d}t,  \qquad  j=1,\ldots,N(n),\vspace*{-2pt}
\]
the estimator of $\mu_k = \P(X(0)\geq u_k)$ based on observations
within $V_{n,j}$, and by\vspace*{-2pt}
\[
\bar{\mu}_{nk} = \frac{1}{N(n)} \sum_{j=1}^{N(n)} \hat{\mu
}^{(j)}_{nk},  \qquad   n\in\mathbb{N},
\ k=1,\ldots,r,\vspace*{-2pt}
\]
the average of these estimators. After all, we define the estimator
$\hat{\Sigma}_n=(\hat{\sigma}_{nlm})_{l,m=1}^r$ for the covariance
matrix $\Sigma$.
Set\vspace*{-2pt}
%
\begin{equation}\label{est:cov3}
\hat{\sigma}_{nlm} = \frac{\nu_d(V_n)}{N(n)-1} \sum_{j=1}^{N(n)}
\bigl(\hat{\mu}^{(j)}_{nl}-\bar{\mu}_{nl}\bigr)\bigl (\hat{\mu}^{(j)}_{nm}-\bar
{\mu}_{nm}\bigr).\vspace*{-2pt}
\end{equation}
We recall the following result.\vspace*{-2pt}
\begin{theorem}[(\cite{Pantle}, Theorem 3)]
Let $X=\{X(t), t\in\R^d\}$ be a strictly stationary random field
such that\vspace*{-2pt}
%
\begin{equation}\label{eq:Cumulant}
\int_{\R^{3d}}\bigl|c_{lm}^{(2,2)}(x,y,z)\bigr|\,\mathrm{d}x\,\mathrm{d}y\,\mathrm{d}z<\infty,\qquad
l,m=1,\ldots, r,\vspace*{-2pt}
\end{equation}
where the fourth-order cumulant function\vspace*{-2pt}
%
\begin{eqnarray*}
c_{lm}^{(2,2)}(x,y,z) &=& \E\bigl([Z_l(0)-\mu_l][Z_m(x)-\mu_m][Z_l(y)-\mu
_l][Z_m(z)-\mu_m]\bigr)
\\[-1pt]
&&{}- \cov_{lm}(x)\cov_{lm}(z-y) - \cov_{ll}(y)\cov_{mm}(x-z) - \cov
_{lm}(z)\cov_{ml}(x-y)\vspace*{-2pt}\vadjust{\goodbreak}
\end{eqnarray*}
and $\cov_{lm}(t)=\cov(\ind\{X(0)\geq u_l\},\ind\{X(t)\geq u_m\})$,
$l,m=1,\ldots,r$. Then $\hat{\Sigma}$ introduced in~\eqref
{est:cov3} is
mean-square consistent.
\end{theorem}

Relation \eqref{eq:Cumulant} holds for a random
field $X$ with finite dependence range. In this case, the estimator
$\hat{\Sigma}_n$ is mean-square consistent.
Among other estimators for the asymptotic covariance matrix, there are
two worth mentioning. One estimator that, under certain assumptions,
meets the conditions of Theorem \ref{th:statvers} is introduced in
\cite{Bulinski2,Bulinski3} and involves local averaging. A major
disadvantage is tedious calculation in the case of a large observation
window. The same problem arises for an estimator based on the
covariance function estimation for the underlying random field; see
\cite{Pantle} (cf. \cite{Ivanov}, Chapter 4).
\vspace*{2pt}\section{Discussion}\vspace*{2pt}\label{discussion}
A very important issue for applications of the estimator $\hat{\Sigma
}_n$ is the choice of an appropriate size of the (e.g., rectangular)
subwindow $V_n$. The subwindow size is related to both the covariance
structure of the considered random field and the size of the
observation window. We will discuss these problems while considering a
simple example. The data used consist of $100$ mutually independent
realizations of stationary and centered Gaussian random field $X$ with
covariance function
\[
\cov(X(0),X(t)) =  \biggl(1-\frac{3 \|t \|_2}{2a}+\frac
{ \|t \|_2^3}{2a^3} \biggr)\ind\{ \|t \|_2\in
[0,a]\},\qquad t\in\R^2,
\]
for some $a>0$ according to the spherical covariance model (see \cite
{Wacker}, page 244), which is often applied in geostatistics. The
correlation range $a$ in our simulation study has to be small enough in
comparison with the size of the observation window to make the CLT
argument work. Here, we take, for example, $a=10$. The fields are
simulated in the observation window $W=[0,2000)\times[0,2000)$ on
the grid with mesh size one. That means every realization
provides $4$ million data points. To generate level sets, we
consider the thresholds $u_1=-1.0$, $u_2=0.0$ and $u_3=1.0$. Then
\[
\Sigma=
\pmatrix{ 4.6432 & &\cr
5.9938 & 10.5564 & \cr
2.7962 & 5.9938 & 4.6432
}
.
\]
An appropriate subwindow size can be found focusing only on the
threshold $u_2=0.0$, since for other threshold values the obtained
results differ from this one only slightly. The estimator provides the
best result for $\Sigma$ as the edge length of the rectangular
subwindow equals $15$. In general, the optimal choice of this length is
an open non-trivial problem.
After this preliminary step, we are able to apply the subwindow
estimator to the simulated data. The following two matrices show
averaged estimation results for $\Sigma$ by means of $\hat{\Sigma}$.
On the
left-hand side, the averaged value of each estimated matrix
element is computed out of 100 samples. On the right-hand side,
the mean error to the theoretical value is provided.
\[
\frac{1}{100}\sum_{k=1}^{100} \hat{\Sigma}_k =
\pmatrix{\displaystyle
4.6556 & & \cr\displaystyle
5.9710 & 10.5524 & \cr\displaystyle
2.8156 & 5.9934 & 4.6762\cr
}
,\qquad
\mathrm{ME} =
\pmatrix{\displaystyle
 0.27\% & & \cr\displaystyle
-0.38\% & -0.04\% & \cr\displaystyle
 0.69\% & \approx0.00\% &  0.71\%\cr
}
.
\]
It would be interesting to propose a statistical hypothesis test based
on the established statistical version of the CLT in order to apply it
to data concerning the paper production. We will deal with this topic
in a separate paper.

\section{Open problems}\label{sectResume}

The research area of limit theorems for level sets of random
surfaces still offers an abundance of open problems. Let us
mention just a few.
It would be desirable to prove limit theorems for joint
distributions of various surface characteristics of different classes of
random fields. For instance, one could consider stable fields. Further
on, one can study random fields possessing
more strong dependence structure; for example, satisfying condition
(A) with $\alpha\leq d$. In this case, the normalizing factors have to
be changed and the limiting distributions can be non-Gaussian. Certain
results for problems of this type can be found in \cite
{Ivanov,Leonenko3}. One could also prove a functional limit theorem for
an innumerable set of thresholds. As our main result could also be
called the CLT for the first Minkowski functional, it might be of
interest to prove limit theorems involving
other Minkowski functionals for level sets such as the boundary length
or the Euler characteristics. It is worth mentioning that, for a
stationary two-dimensional Gaussian field, this has already been done
for the second Minkowski functional in \cite{kratzleon01}.

\section*{Acknowledgements}
The authors are grateful to Professor Richard A. Davis, the Associate
Editor and the referees for valuable remarks and suggestions permitting
them to improve the exposition of the paper.
Alexander Bulinski's work was partially supported by RFBR Grant
10-01-00397-a.

%

\printhistory

\end{document}